\documentclass[aos,seceqn,nameyear,dvips]{arximspdf}
\usepackage{dcolumn}
\usepackage{stfloats}
\usepackage{graphicx}

\doi{10.1214/10-AOS794}
\volume{38}
\issue{5}
\pubyear{2010}
\firstpage{2652}
\lastpage{2677}

\makeatletter

\fnbelowfloat

\newcolumntype{d}[1]{D{.}{.}{#1}}

\newtheorem{theorem}{Theorem}
\newtheorem{lemma}{Lemma}

\newproclaim{Example}{Example}

\makeatother

\begin{document}
\begin{frontmatter}

\title{Modeling the variability of rankings}
\runtitle{Ranking}

\begin{aug}
\author[A]{\fnms{Peter} \snm{Hall}\ead[label=e1]{halpstat@ms.unimelb.edu.au}} and
\author[A]{\fnms{Hugh} \snm{Miller}\corref{}\ead[label=e2]{h.miller@ms.unimelb.edu.au}}
\runauthor{P. Hall and H. Miller}
\affiliation{University of Melbourne}
\address[A]{Department of Mathematics and Statistics\\
University of Melbourne\\
Melbourne, VIC 3010\\
Australia\\
\printead{e1}\\
\phantom{E-mail: }\printead*{e2}} 
\end{aug}

\received{\smonth{6} \syear{2009}}
\revised{\smonth{10} \syear{2009}}

%
\begin{abstract}
For better or for worse, rankings of institutions, such as
universities, schools and hospitals, play an important role today in
conveying information about relative performance. They inform policy
decisions and budgets, and are often reported in the media. While
overall rankings can vary markedly over relatively short time periods,
it is not unusual to find that the ranks of a small number of ``highly
performing'' institutions remain fixed, even when the data on which the
rankings are based are extensively revised, and even when a large
number of new institutions are added to the competition. In the
present paper, we endeavor to model this phenomenon. In particular,
we interpret as a random variable the value of the attribute on which
the ranking should ideally be based. More precisely, if $p$ items are
to be ranked then the true, but unobserved, attributes are taken to be
values of $p$ independent and identically distributed variates.
However, each attribute value is observed only with noise, and via a
sample of size roughly equal to $n$, say. These noisy approximations to
the true attributes are the quantities that are actually ranked. We
show that, if the distribution of the true attributes is light-tailed
(e.g., normal or exponential) then the number of institutions
whose ranking is correct, even after recalculation using new data and
even after many new institutions are added, is essentially fixed.
Formally, $p$ is taken to be of order $n^C$ for any fixed $C>0$, and
the number of institutions whose ranking is reliable depends very
little on $p$. On the other hand, cases where the number of reliable
rankings increases significantly when new institutions are added are
those for which the distribution of the true attributes is relatively
heavy-tailed, for example, with tails that decay like $x^{-\alpha}$ for
some $\alpha>0$. These properties and others are explored analytically,
under general conditions. A numerical study links the results to
outcomes for real-data problems.
\end{abstract}

%

\end{frontmatter}

\section{Introduction}\label{sec1}

There are many contemporary settings in which ranking plays an
important role. For example, universities, schools and hospitals are
regularly ranked in a variety of contexts, the results of which
typically generate interest and can often drive policy decisions. In
many of these situations, a given ranking can carry a high degree of
uncertainty, with this effect particularly pronounced in high-dimensional
cases; that is, where there are very many populations or
institutions to be ranked.

Despite this, one feature of many rankings reported over time is that
the ordering at the extreme top or bottom remains relatively invariant.
For example, in the THE-QS university
rankings,\footnote{\href{http://www.topuniversities.com}{www.topuniversities.com}.}
Harvard University has ranked first for each
of the years 2005--2008, while New York University's rankings are 56,
43, 49 and 40. If we believe that the observed data used for ranking
are measures of true underlying values, distorted by noise, then we can
reinterpret this behavior as a tendency to obtain correct rankings at
extremes, but not otherwise. It is this phenomenon that we explore in
this paper, using both theoretical and numerical arguments.

Intuitively, this behavior has a natural explanation. Those scores at
the extreme of a range are more likely to be sufficiently ``spaced
out'' to overcome the problems of data noise, whereas less extreme
scores are likely to be bunched more closely together. We introduce
models that describe this behavior and explore their properties.
Related to this, it turns out that one important consideration for
correct ranking at the extremes is whether the possible scores used for
ranking have infinite support but nevertheless have light tails. If
this is the case and the tail of the distribution of the underlying
scores is smooth, we can expect accurate ranking of the top portion of
the institutions, even when dimension is very large. Moreover, even
when the support is bounded, there remains potential for correct
ranking at extremes, although now there is greater likelihood that the
ranking will change if new institutions are added. Such results have a
variety of practical implications; we briefly present two of these
here, with more detail provided in the numerical section.
\begin{Example}[(University rankings)]\label{Example1}
Suppose we attempt to rank
universities and other research institutions by counting how many
papers their faculty members publish in
\textit{Nature}\footnote{\href{http://www.nature.com/nature/index.html}{www.nature.com/nature/index.html}.} each
year. This is a high-dimensional
example due to the large number of institutions competing
to be published. Figure \ref{fig1} shows the ranking of the top 50
institutions on this measure. The institutions are aligned along the
horizontal axis, with the each dot denoting the point estimate of the
rank and the vertical line a corresponding estimated 90\% prediction
interval. The four plots show how the confidence intervals change as we
increase the number of years, $n$, of data used for the ranking.

%
\begin{figure}

\includegraphics{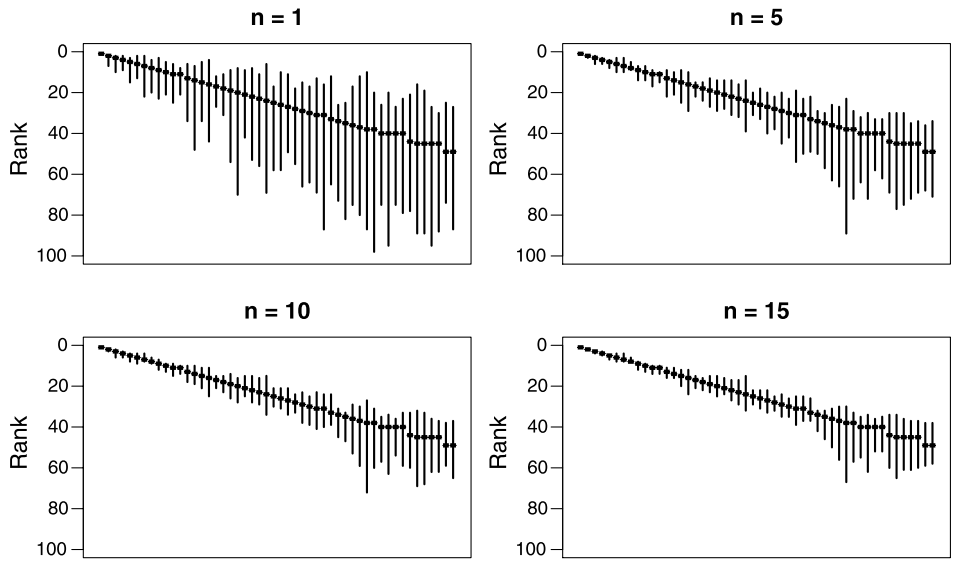}

\caption{Prediction intervals for top-ranked universities
based on publications in Nature, averaged over various numbers of
years.}\label{fig1}
\end{figure}

The two main observations are that the prediction intervals are widest
when a smaller number of years are considered, and that the prediction
intervals for the highest ranked universities are the smallest. In
fact, the intervals are small enough in the extremes to give us genuine
confidence in that aspect of the ranking. Even when $n=1$, we can be
reasonably sure that the top ranked institution (Harvard University) is
in fact ranked correctly. When $n=15$, the top four universities are
known with a high degree of certainty, and the next set of ten or so is
fairly stable too. Thus, it is possible to have correctness in the
upper extreme of this ranking, even when the lower ranks remain highly
variable. In the present paper, we model this phenomenon by addressing
the underlying stochastic properties of the institutions; the data
provide only a noisy measure of this random process, and we assess the
impact of the noise on the ranking.
\end{Example}
\begin{Example}[(Microarray data)]\label{Example2}
We take the colon microarray data first analysed by
\citet{Alonetal99}. It consists of 62 observations in total, each
of which indicates either a normal colon or a tumor. For each
observation, there are also expression levels for $p=2000$ genes. It is
of interest to determine which genes are most closely related to the
response, so that they can be investigated further. This of course
amounts to a ranking and we are interested in stability at the extreme,
since we seek only a small number of genes. Here, the genes are ranked
based on the Mann--Whitney U test statistic, which is a nonparametric
assessment of the difference between the two distributions.

%
\begin{figure}

\includegraphics{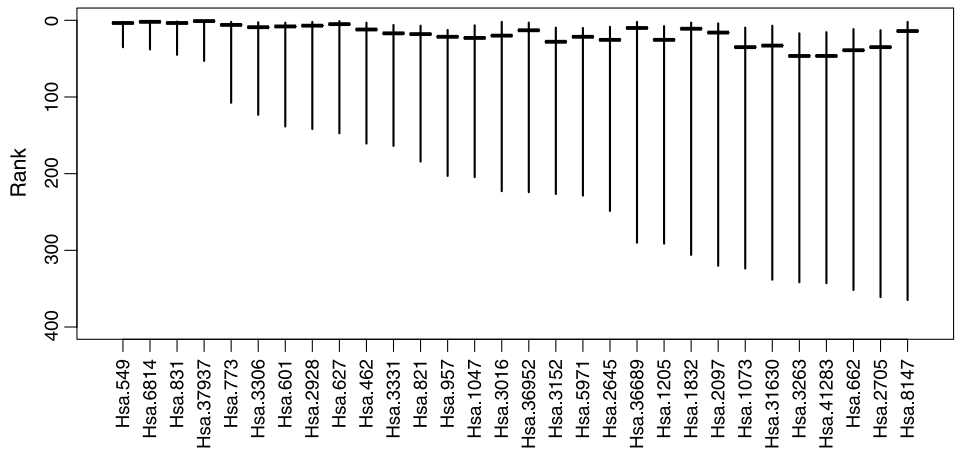}

\caption{Prediction intervals for top-ranked genes in colon dataset.}
\label{fig2}
\end{figure}

Figure \ref{fig2} plots the top 30 genes, ranked by the lower tail of an
estimated 90\% prediction interval, rather than the point estimate of
the rank. In this situation, we cannot authoritatively conclude that
any of the top genes are ranked exactly correctly, but the top four
genes appear much more stable than the others. This stability is highly
important; if the length of all prediction intervals were roughly the
same as the average length (1400 genes), then there would be little
hope of discovering useful genes from such datasets.

There is a literature on the bootstrap in connection with rankings. See
\citet{GS96}, who discuss bootstrap methods for
constructing prediction intervals for rankings;
\citet{LL96}, who address bootstrap methods for ranking the
performance of
doctors; \citet{CB03}, \citet{HMZ05} and \citet{TB05},
who take up the problem of bootstrap methods for
ranked set sampling; \citet{Mukherjeeetal03}, who develop methods for
gene ranking using bootstrapped $p$-values; and \citet{XSZ09} and
\citet{HM10}, who focus on consistent bootstrap methods for
assessing rankings. More generally, there is a vast literature on
ranking problems in statistics, and we cite here only the more relevant
items since 2000. Joe (\citeyear{J00}, \citeyear{J01}) discusses
ranking problems in
connection with random utility models, and points to connections to
multivariate extreme value theory. \citet{MM03} develop
mixture-based models for rankings. \citet{M03} and \citet
{Barkeretal05}
treat methods for ranking football players. \citet{MS05}
study the problem of ranking immunisation coverage in U.S. states. Brijs
et al. (\citeyear{Brijsetal06}, \citeyear{Brijsetal07}) introduce
Bayesian models for the ranking of
hazardous road sites, with the aim of better scheduling road safety
policies. \citet{CSW06} discuss ranking accuracy in ranked-set
sampling methods, and \citet{OrS07} examine the
accuracy of gene rankings in high-dimensional problems involving
genomic data. \citet{N06} addresses the reliability of performance
rankings. \citet{CS07} and \citet{QBL07} discuss
ways of constructing rankings.

Section \ref{sec2} describes our model for the ranking problem, and discusses
the main properties of this framework. The formal theoretical results
which underpin the discussion in Section \ref{sec2} are given in
Section \ref{sec3}.
Section \ref{sec4} presents simulated and real-data numerical work, including
details on the examples presented above. Technical proofs are deferred
to Section \ref{sec5}.
\end{Example}

\section{Model}\label{sec2}

We consider a set of underlying parameters $\theta_1,\ldots,\theta
_p$ corresponding to the objects to be ranked, hereafter referred to as
items. The error in the estimation is controlled by the number of
observed data points, $n$. In our analysis, we take $p=p(n)$ to diverge
with $n$ as the latter increases. An obvious difficulty here is in
establishing where the newly added items should fit into the ranking.
A~natural solution is to take the $\theta_j$'s to be randomly generated
from some distribution function. In the setup below, we interpret the
$\Theta_j$'s as values of means; see the end of this section for
generalizations.

Let $\Theta_1,\ldots,\Theta_p$ denote independent and identically
distributed random variables, and write
%
%
\begin{equation}\label{eq1.1}
\Theta_{(1)}\leq\cdots\leq\Theta_{(p)}
\end{equation}
for their ordered values. There exists a permutation $R=(R_1,\ldots
,R_p)$ of $(1,\ldots,p)$ such that $\Theta_{(j)}=\Theta_{R_j}$ for
$1\leq
j\leq p$. If the common distribution of the $\Theta_j$'s is continuous,
then the inequalities in (\ref{eq1.1}) are all strict and the
permutation is unique.

We typically do not observe the $\Theta_j$'s directly, only in terms of
noisy approximations which can be modelled as follows. Let
$Q_i=(Q_{i1},\ldots,Q_{ip})$ denote independent and identically
distributed random $p$-vectors with finite variance and zero mean,
independent also of $\Theta=(\Theta_1,\ldots,\Theta_p)$. Suppose
we observe
%
%
\begin{equation} \label{eq1.2}
X_i=(X_{i1},\ldots,X_{ip})
=Q_i+\Theta
\end{equation}
for $1\leq i\leq n$. The mean vector
%
%
\begin{equation} \label{eq1.3}
{\bar X}=({\bar X}_1,\ldots,{\bar X}_p)
=\frac{1}{n}\sum_{i=1}^n X_i={\bar Q}+\Theta
\end{equation}
is an empirical approximation to $\Theta$. (Here, ${\bar
Q}=n^{-1}\sum_i Q_i$
equals the mean of the $p$-vectors $Q_i$.) The components of ${\bar X}$ can
also be ranked, as
%
%
\begin{equation}\label{eq1.4}
{\bar X}_{(1)}\leq\cdots\leq{\bar X}_{(p)},
\end{equation}
and there\vspace*{1pt} is a permutation ${\widehat R}_1,\ldots,{\widehat R}_p$ of
$1,\ldots,p$
such that ${\bar X}_{(j)}=X_{{\widehat R}_j}$ for each~$j$. If the
common distribution
of the $\Theta_j$'s is continuous then, regardless of the distribution of
the components of $Q_i$, the inequalities in (\ref{eq1.4}) are strict
with probability~1.

The permutation ${\widehat R}=({\widehat R}_1,\ldots,{\widehat R}_p)$
serves as an
approximation to $R$, and we wish to determine the accuracy of that
approximation. In particular, for what values of $j_0=j_0(n,p)$, and
for what relationships between $n$ and $p$, is it true that
%
%
\begin{equation} \label{eq1.5}
P ({\widehat R}_j=R_j \mbox{ for } 1\leq j\leq j_0 )\to1
\end{equation}
as $n$ and $p$ diverge? That is, how deeply into the ranking can we go
before the connection between the true ranking and its empirical form
is seriously degraded by noise?

The answer to this question depends to some degree on the extent of
dependence among the components of each $Q_i$. To elucidate this point,
let us consider the case where all the components of $Q_i$ are
identical; this is an extreme case of strong dependence. Then the
components of ${\bar Q}$ are also identical. Clearly, in this setting
${\widehat R}
_j=R_j$ for each $j$, and so (\ref{eq1.5}) holds in a trivial and
degenerate fashion. Other strongly dependent cases, although not as
clear-cut as this one, can also be shown to be ones where ${\widehat R}_j=R_j$
with high probability for many values of $j$.

The case which is most difficult, that is, where the strongest
conditions are needed to ensure that (\ref{eq1.5}) holds, occurs when
the components of $Q_i$ are independent. To emphasize this point we
give sufficient conditions for (\ref{eq1.5}), and show that when the
components of each $Q_i$ are independent, those conditions are also
necessary. Our arguments can be modified to show that the conditions
continue to be necessary under sufficiently weak dependence, for
example if the components are $m$-dependent where $m=m(n)$ diverges
sufficiently slowly as $n$ increases.

The assumptions under which (\ref{eq1.5}) holds are determined mainly
by the lower tail of the common distribution of the $\Theta_j$'s. If that
distribution has an exponentially light left-hand tail, for example, if
the tail is like that of a normal distribution, then a sufficient
condition for (\ref{eq1.5}) is that $j_0$ should increase at a
strictly slower rate than $n^{1/4} (\log n)^c$, where the constant
$c$, which can be either positive or negative, depends on the rate of
decay of the exponential lower tail of the distribution of $\Theta$. For
example, $c=0$ if the distribution decays like $e^{-|x|}$ in the lower
tail, and $c=-{\frac{1}{4}}$ if it is normal. As indicated in the previous
paragraph, the condition $j_0 =o\{n^{1/4} (\log n)^c\}$ is also
necessary for (\ref{eq1.5}) if the components of the $Q_i$'s are independent.

These results have several interesting aspects, including: (a) the
exponent ${\frac{1}{4}}$ in the condition $j_0=o\{n^{1/4} (\log n)^c\}
$ does
not change among different types of distribution with exponential
tails; (b) the exponent is quite small, implying that the empirical
rankings ${\widehat R}_j$ quite quickly become unreliable as predictors
of the
true rankings $R_j$; and (c) the critical condition $j_0=o\{n^{1/4}
(\log n)^c\}$ does not depend on the value of $p$. (We assume that $p$
diverges at no faster than a polynomial rate in~$n$, but we impose no
upper bound on the degree of that polynomial.)

The condition on $j_0$ such that (\ref{eq1.5}) holds changes in
important ways if the lower tail of the distribution of the $\Theta_j$'s
decays relatively slowly, for example, at the polynomial rate
$x^{-\alpha}$
as $x\to\infty$. Examples of this type include Pareto, nonnormal stable
and Student's $t$ distributions, and more generally, distributions with
regularly varying tails. Here a sufficient condition for (\ref{eq1.5})
to hold is $j_0=o\{(n^{\alpha/2} p)^{1/(2\alpha+1)}\}$, and this assumption
is necessary if the components of the $Q_i$'s are independent. In this
setting, unlike the exponential case, the value of dimension, $p$,
plays a major role in addition to the sample size, $n$, in determining
the number of reliable rankings.

In practical terms, a major way in which this heavy-tailed case differs
from the light-tailed setting considered earlier is that if a
polynomially large number of new items are added to the competition in
the heavy-tailed case, and all items are reranked, the results will
change significantly and the number of correct rankings will also alter
substantially. By way of contrast, if a polynomially large number of
new items are added in the light-tailed, or exponential, case then
there will again be many changes to the rankings, but now there will be
relatively few changes to the number of items that are correctly ranked.

The exponential case can be regarded as the limit, as $\alpha\to
\infty$, of
the polynomial case. More generally, note that as the left-hand tail of
the common distribution of the $\Theta_j$'s becomes heavier, the value of
$j_0$ can be larger before (\ref{eq1.5}) fails. That is, if the
distribution of the $\Theta_j$'s has a heavier left-hand tail then the
empirical rankings ${\widehat R}_j$ approximate the true rankings $R_j$
for a
greater number of values of $j$, before they degenerate into noise.

The analysis above has focused on cases where the ranks of the $X_j$'s
are estimated by ranking empirical means of noisy observations of those
quantities; see~(\ref{eq1.4}). However, similar results are obtained
if we rank other measures of location. Such a measure need only satisfy
moderate deviation properties similar to (\ref{eqP.5}) and (\ref
{eqP.6}) in the proof of Theorem \ref{theo1}. Thus, the results are
applicable to
a wide range of ranking contexts. For example, $L_q$ location
estimators for general $q\geq1$ enjoy moderate deviation properties
under appropriate assumptions. Therefore, if we take the variables
$Q_{ij}$ to have zero median, rather than zero mean, and continue to
define $X_i$ by (\ref{eq1.2}) but replace the ranking in (\ref
{eq1.4}) by a ranking of medians, then the results above and those in
Section \ref{sec3} continue to hold, modulo changes to the regularity
conditions. Other suitable measures include the Mann--Whitney test used
in the genomic example, quantiles and some correlation-based measures.

The model suggested by (\ref{eq1.2}), where data on $\Theta$ arise in
the form of $p$-vectors $X_1,\ldots,X_n$, is attractive in a number of
high-dimensional settings, for example, genomics. There, the $j$th
component $X_{ij}$ of $X_i$ would typically represent the expression
level of the $j$th gene of the $i$th individual in a sample. However,
in other cases the means ${\bar X}_1,\ldots,{\bar X}_p$ at (\ref
{eq1.3}), or
medians or other location estimators, might be computed from quite
different datasets, one for each component index $j$. Moreover, those
datasets might be of different sizes, $n_1,\ldots,n_p$ say, and then
the argument that they arise naturally in the form of vectors would be
inappropriate. This can happen when data are used to rank items, for
example schools where the ranking is based on individual student
performance. The conclusions discussed earlier in this section, and the
theoretical properties developed in Section \ref{sec3} below, continue
to apply
in this case provided there is an ``average'' value, $n$ say, of the
$n_j$'s which represents all of them, in the sense that
%
%
\begin{equation}\label{eq1.6}
n=O \Bigl(\min_{1\leq j\leq p} n_j \Bigr) \quad\mbox{and}\quad
\max_{1\leq j\leq p} n_j=O(n)
\end{equation}
as $n$ diverges. Additionally, in such cases it is often realistic to
make the assumption that the corresponding centred means (or medians,
etc.) ${\bar Q}_j=n^{-1}\sum_i Q_{ij}$ are stochastically independent
of one
another, and so the particular results that are valid in this case are
immediately available.

The distribution of the $\Theta_j$'s has been taken to be continuous. This
is usually appropriate although there can be contexts in which the
distribution is discrete. Note that assumption of discreteness of the
$\Theta_j$'s is different from that of discreteness of the
observations $X_{ij}$. In such cases, the analysis still holds, except
that allowance must be made for ties (any reordering of tied $\Theta_j$'s
is still ``correct''), and the tail density assumptions should be
characterized in integral form.

The model has been set up so that it focuses on the populations with
lowest parameters $\Theta_j$. Obviously, similar arguments apply to the
largest parameters too, so the results are applicable to both the most
highly and lowly ranked populations.

\section{Theoretical properties}\label{sec3}

For the most part, we shall assume one of two types of lower tail for
the common distribution function, $F$, of the random variables $\Theta
_j$: either it decreases exponentially fast, in which case we suppose
that $F(-x)\asymp x^\beta\exp(-C_0 x^\alpha)$ as $x\to\infty$, where
$\alpha>0$ and $-\infty<\beta<\infty$; or it decreases
polynomially fast,
in which case $F(-x)\asymp x^{-\alpha}$ as $x\to\infty$, where
$C_0,\alpha>0$.
[The notation $f(x)\asymp g(x)$, for positive functions $f$ and $g$,
will be taken to mean that $f(x)/g(x)$ is bounded away from zero and
infinity as $x\to\infty$.] The former case covers distributions such
as the
normal, exponential and Subbotin; the latter, distributions such as the
Pareto, Student's $t$ and nonnormal stable laws (e.g., the Cauchy).

It is convenient to impose the shape constraints on the densities,
which we assume to exist in the lower tail, rather than on the
distribution functions. Therefore, we assume that one of the following
two conditions hold as $x\to\infty$:
%
%
\begin{eqnarray}
\label{eq2.1}
(d/dx) F(-x)&\asymp& (d/dx) x^\beta\exp(-C_0 x^\alpha) ,\\
\label{eq2.2}
(d/dx) F(-x)&\asymp& (d/dx) x^{-\alpha} .
\end{eqnarray}
In both (\ref{eq2.1}) and (\ref{eq2.2}), $\alpha$ must be strictly
positive, but $\beta$ in (\ref{eq2.1}) can be any real number. The
constant $C_0$ in (\ref{eq2.1}) must be positive. We assume too that
%
%
\begin{equation}\label{eq2.3}
\begin{tabular}{p{11cm}}
for fixed constants $C_1,\ldots,C_5>0$, where $C_2>2 (C_1+1)$ and
$C_4<C_5$, $p=O(n^{C_1})$ as $n\to\infty$, and, for each $j\geq1$,
$E|Q_j|^{C_2}\leq C_3$, $E(Q_j)=0$, and
$E(Q_j^2)\in[C_4,C_5]$.
\end{tabular}\hspace*{-28pt}
\end{equation}

Recall from Section \ref{sec1} that we wish to examine the probability
that the
true ranks~$R_j$, and their estimators ${\widehat R}_j$, are identical
over the
range $1\leq j\leq j_0$. We consider both $j_0$ and $p$ to be functions
of $n$, so that the main dependent variable can be considered to
be $n$. With this interpretation, define
%
%
\begin{eqnarray} \label{eq2.4}
\nu_{\exp}&=&\nu_{\exp}(n)=n^{1/4} (\log n)^{\{(1/\alpha)-1\}/2}
,\nonumber\\[-8pt]\\[-8pt]
\nu_{\mathrm{pol}}&=&\nu_{\mathrm{pol}}(n)=(n^{\alpha/2}
p)^{1/(2\alpha+1)} ,\nonumber
\end{eqnarray}
where the subscripts denote ``exponential'' and ``polynomial,''
respectively, and refer to the respective cases represented by (\ref
{eq2.1}) and (\ref{eq2.2}). In the theorem below, we impose the
additional condition that, for some $\varepsilon>0$,
%
%
\begin{equation} \label{eq2.4b}
n=O(p^{4-\varepsilon}).
\end{equation}
This restricts our attention to problems that are genuinely
high dimensional, in the sense that, with probability converging to 1,
not all the rankings are correct. Cases where $p$ diverges sufficiently
slowly as a function of $n$ are easier and will generally permit all
ranks to be correctly determined with high probability. Assumption
(3.5) is also very close, in both the exponential and polynomial cases,
to the basic condition $j_0\leq p$, as can be seen via a little
analysis starting from (3.6) and (3.7) in the respective cases; yet, at
the same time, (3.5) is suitable to both cases, and so helps to unify
our account of their properties. Note too that (3.5) implies that, in
both the exponential and polynomial cases, $\nu_{\exp}=O(p^{1-\delta
})$ and $\nu_{\mathrm{pol}}=O(p^{1-\delta})$ for some $\delta>0$.
\begin{theorem}\label{theo1}
Assume (\ref{eq2.3}), (\ref{eq2.4b}) and that either
\textup{(a)} (\ref{eq2.1}) or \textup{(b)} (\ref{eq2.2}) holds.
In case \textup{(a)}, if
%
%
\begin{equation} \label{eq2.5}
j_0=o(\nu_{\exp})
\end{equation}
as $n\to\infty$ then (\ref{eq1.5}) holds. Conversely, when the components
of the vectors $Q_i$ are independent, (\ref{eq2.5}) is necessary
for (\ref{eq1.5}). In case \textup{(b)}, if
%
%
\begin{equation} \label{eq2.6}
j_0=o(\nu_{\mathrm{pol}}) ,
\end{equation}
then (\ref{eq1.5}) holds. Conversely, when the components of the
vectors $Q_i$ are independent, (\ref{eq2.6}) is necessary for (\ref{eq1.5}).
\end{theorem}

It can be deduced from Theorem \ref{theo1} that when a new item (e.g., an
institution) enters the competition that leads to the ranking, we are
still able to rank the top $j_0$ institutions correctly. In this sense,
the institutions that make up the cohort of size $j_0$ do not need to
be fixed.

It is also of interest to consider cases where the common distribution,
$F$, of the $\Theta_j$'s is bounded to the left, for example, where
$F(x)\asymp x^\alpha$ as $x\downarrow0$. However, it can be shown
that in
this context, unless $p$ is constrained to be a sufficiently low degree
polynomial function of $n$, very few of the estimated ranks ${\widehat R}_j$
will agree with the correct values $R_j$.

To indicate why, we first recall the model introduced in Section \ref{sec1},
where the estimated ranks ${\widehat R}_j$ are derived by ordering the values
of ${\bar Q}_j+\Theta_j$. Here ${\bar Q}_j=n^{-1}\sum_{1\leq i\leq
n} Q_{ij}$
is the average value of $n$ independent and identically distributed
random variables with zero mean. Therefore the means, ${\bar Q}_j$, are of
order $n^{-1/2}$. By way of contrast, if we take $\alpha=1$ in the formula
$F(x)\asymp x^\alpha$ as $x\downarrow0$, for example, if $F$ is the
uniform distribution on $[0,1]$, then the spacings of the order
statistics $\Theta_{(1)}\leq\cdots\leq\Theta_{(p)}$ are
approximately of
size $p^{-1}$. (More concisely, they are of size $Z/p$ where $Z$ has an
exponential distribution; an independent version of $Z$ is used for
each spacing.) Therefore, if $p$ is of larger order than $n^{1/2}$ then
the errors of the ``estimators'' ${\bar Q}_j+\Theta_j$ of $\Theta
_j$, for
$1\leq j\leq p$, are an order of magnitude larger than the spacings
among the $\Theta_j$'s. This can make it very difficult to estimate the
ranks of the $\Theta_j$'s from the ranks of values of ${\bar
Q}_j+\Theta_j$.
Indeed, it can be shown that, in the difficult case where the
components of the $Q_i$'s are independent, and even for fixed $j_0$, if
$\alpha=1$ and $p$ is of larger order than $n^{1/2}$ then in contrast
to (\ref{eq1.5}),
%
%
\begin{equation} \label{eq2.7}
P ({\widehat R}_j=R_j \mbox{ for } 1\leq j\leq j_0 )\to0 .
\end{equation}

This explains why, when $F(x)\asymp x^\alpha$, it can be quite rare for
the estimated ranks ${\widehat R}_j$ to match their true values.
Indeed, no
matter what the value of $\alpha$ and no matter what the value of $j_0$,
property (\ref{eq1.5}) will typically fail to hold unless $p$ is no
greater than a sufficiently small power of $n$, in particular unless
$p=o(n^{\alpha/2})$, as the next result indicates. Thus, the differences
between the cases of bounded and unbounded distributions are stark, as
can be seen by contrasting Theorem \ref{theo1} with the properties
described below.
\begin{theorem}\label{theo2}
Assume that $(d/dx) F(x)\asymp x^{\alpha-1}$ as
$x\downarrow0$, where $\alpha>0$, and that (\ref{eq2.3}) holds.
Part \textup{(a)}: instances where (\ref{eq1.5}) holds and $p^2/n^\alpha\to0$.
Under the latter condition, \textup{(i)} if $\alpha<{\frac{1}{2}}$ then (\ref{eq1.5})
holds even for $j_0=p$; \textup{(ii)} if $\alpha={\frac{1}{2}}$ then
(\ref{eq1.5}) holds
provided that
%
%
\begin{equation}\label{eq2.8}
(\log j_0)^{2\alpha} (p^2/n^\alpha)\to0 ;
\end{equation}
and \textup{(iii)} if $\alpha>{\frac{1}{2}}$ then (\ref{eq1.5}) holds
provided that
%
%
\begin{equation}\label{eq2.9}
j_0=o \bigl\{(n^{\alpha/2}/p)^{1/(2\alpha-1)} \bigr\}.
\end{equation}
Part \textup{(b)}: converses to \textup{(a)(ii)} and \textup{(a)(iii)}.
If $p^2/n^\alpha\to0$ and the components of the vectors $Q_i$
are independent then, if (\ref{eq1.5}) holds,
so too does (\ref{eq2.8}) (if $\alpha={\frac
{1}{2}}$) or (\ref
{eq2.9}) (if $\alpha>{\frac{1}{2}}$). Part \textup{(c)}: instances where
(\ref{eq2.7})
holds. If $\alpha>0$ and $p^2/n^\alpha\to\infty$, and if the
components of the
vectors $Q_i$ are independent, then (\ref{eq2.7}) holds even for $j_0=1$.
\end{theorem}

The proof of Theorem \ref{theo2} is similar to that of Theorem \ref
{theo1}, and so is
omitted. Theorem \ref{theo1} is derived in Section \ref{sec5}. Both
results continue to
hold if the sample from which ${\bar X}_j$ is computed is of size $n_j$ for
$1\leq j\leq p$, rather than $n$, provided that (\ref{eq1.6}) holds.

\section{Numerical properties}\label{sec4}

This section discusses three real-data and three simulated examples
linked to the theoretical properties in Section \ref{sec3}. The real-data
examples make use of the bootstrap to create prediction intervals
[\citet{XSZ09}, \citet{HM10}]. In each simulated example, the
error is relatively light-tailed, and any discussion of tails refers to
the distribution of the $\Theta_j$'s. In our real-data examples, the noise
has been averaged and so it is also generally light-tailed. Thus, any
heavy-tailed behavior present in the real-data examples is likely to be
due to heavy tails of the distribution of the $\Theta_j$'s, rather than
the noise.
\setcounter{Example}{0}
\begin{Example}[(Continued)]
The originating institutions of
\textit{Nature} articles were obtained using the ISI Web of knowledge
database\footnote{\href{http://www.isiknowledge.com}{www.isiknowledge.com}.} for each of the years 1999
through 2008. A point ranking was obtained by taking the average number
of articles published per year. Of course, there are implicit
simplifying assumptions in doing this, most significantly concerning
the independence of articles between years, and the stationarity of
means time. These assumptions appear reasonable in context, and are
consistent with most publication-based analyses.

When constructing prediction intervals the bootstrap resamples for each
institution were drawn independently, conditional on the data. [See
\citet{HM10}.] The number of observations in the resample can
be varied to create different time windows, as illustrated in
Figure \ref{fig1}. The most natural question from a ranking
correctness viewpoint is determining the behavior at the right tail;
there are many institutions with mean at or near the hard threshold of
zero, so there is little hope for ranking correctness in the left tail.
Furthermore, the right tail appears to be long. Harvard University has
an average of 67.5 papers per year, followed by means of 34.6, 29.6
and 28.2 for Berkeley, Stanford and Cambridge, respectively.


%
\begin{figure}[b]

\includegraphics{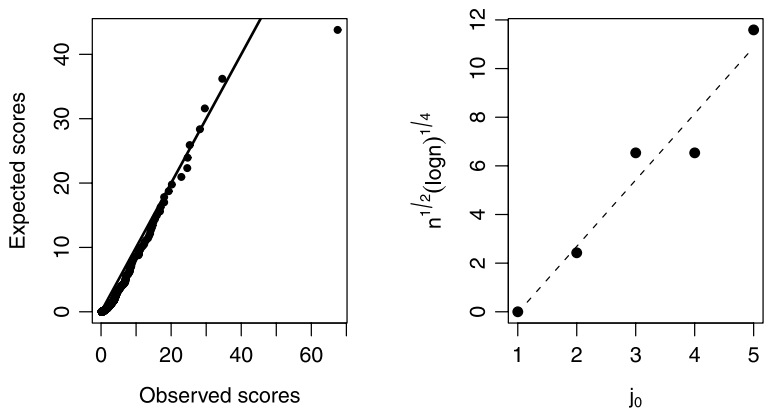}

\caption{The left panel is a QQ plot for the Nature data
against the exponential distribution. The right panel plots a transform
of the number of years of data required to rank $j_0$ institutions
correctly for various $j_0$.}\label{fig3b}
\end{figure}

A natural question to ask is what the tail shape for this example might
be. Approaches to estimating the shape parameter of a distribution with
regularly varying tails, such as the method of \citet{H75}, are
unstable for these data; the number of extreme data for which a linear
fit is plausible is very small, implying that the decay rate is faster
than polynomial. Indeed, the left panel of Figure \ref{fig3b} shows
the QQ plot of the observed data against a random variable with
distribution function $F(x) = 1 - \exp(-0.85 x^{1/2})$, which
suggests that an exponential tail might be reasonable for the data. If
this is the case then the number of institutions that we expect to be
ranked correctly should depend, to first order, only on $n$, not on
$p$, and be of order up to $n^{1/4} (\log n)^{1/2}$. One way to
explore this further is to take $j_0$ as given, and to resample from
the data, seeking, for example, the number of years, $n$, needed to
obtain correct ranking of the first $j_0$ institutions at least 90\% of
the time. A plot of $j_0$ against $n^{1/4} (\log n)^{1/2}$ should be
roughly linear. The right-hand panel of Figure \ref{fig3b} plots
results of this experiment and appears to support the hypothesis. The
flatness between $j_0=3$ and $j_0=4$ indicates that these two
institutions are quite difficult to separate from each other.
\end{Example}
\begin{Example}[(Continued)]
The Mann--Whitney test statistic can be written as
\[
\max\biggl\{ \sum_{i,j} I(x_i < y_j), \sum_{i,j} I(x_i > y_j) \biggr\},
\]
where the $x_i$'s and $y_j$'s are the observed values of the two samples.
Notice that this statistic will have a hard lower threshold at
$n_1n_2/2$, where $n_1$ and $n_2$ are the sizes of the two classes.
Here, like the previous example, when the distributions differ only in
location the difference has to be quite large to be detectable.
Figure \ref{fig4} shows the estimated density as well as the truncated
normal density, which is the distribution that the scores would have if
none of the genes had systematically different means for the two
classes. This suggests that an assumption that the majority of genes is
unrelated to whether the tissue is tumorous is not valid here.

%
\begin{figure}[b]

\includegraphics{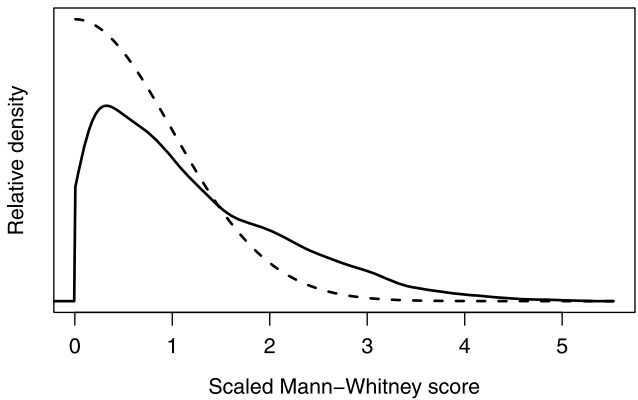}

\caption{Estimated sampling density genes under the
Mann--Whitney test for colon data.}\label{fig4}
\end{figure}

Bootstrapped versions of the dataset with different choices for $n$
were created to indicate how many observations we need to obtain
reasonable confidence in a ranking. Table \ref{tabcolon} shows the
probability that the set of the top $j$ genes is identified correctly
out of the 2000 for various $j$ and $n$. Note that this is a slightly
different statistic from the one in (\ref{eq1.5}), since we allow any
permutation of the top $j$ genes to be detected. The results suggest
that we have nearly a 50\% chance of detecting the top gene if $n=250$,
and a 20\% chance of correctly choosing the top four. The upper tail
for this dataset again appears relatively light; the model $F(x)=1-\exp
\{-0.19 (x-1)^2\}$, for $x>1$, produces a good fit to the upper tail.

%
\begin{table}
\caption{Probability that set of top $j$ genes is correct for
colon data}\label{tabcolon}
\begin{tabular*}{\tablewidth}{@{\extracolsep{\fill}}lccccc@{}}
\hline
&\multicolumn{5}{c@{}}{$\bolds n$} \\ [-4pt]
&\multicolumn{5}{c@{}}{\hrulefill} \\
\multicolumn{1}{@{}l}{$\bolds j$} &
\textbf{62} & \textbf{100} & \textbf{150}
& \textbf{200} & \textbf{250} \\ \hline
$1$ &0.251&0.326&0.437&0.446&0.490\\
$2$ &0.067&0.109&0.166&0.218&0.277\\
$4$ &0.022&0.054&0.094&0.163&0.193\\
$6$ &0.007&0.018&0.035&0.040&0.068\\
\hline
\end{tabular*}
\end{table}

Theorem \ref{theo1} suggests that these probabilities should not depend
on the
choice of $p$. We can obtain a sense of this by randomly sampling,
without replacement, $p=500$ or $p=1000$ genes from the original
$p=2000$, for each simulation; and recalculating the values in
Table \ref{tabcolon}. For $j=4$ and $n=250$, the respective probabilities were 0.183
and 0.170, quite close to the value 0.193 observed for $p=2000$.
While the equivalence appears good for $j \ge4$, there are larger
departures for $j = 1$ or $2$, where the initial results for this
particular realization tend to distort the calculation.
\end{Example}
\begin{Example}[(School rankings)]\label{Example3}
A third example of accuracy in the extremes of a ranking is based on
student performance at 75 private schools in NSW, Australia. For each
school the number of final year exams taken, and the number of these
where a score of at least 90\% was achieved, were recorded. The
proportion of exams where 90\% or more was scored can be used to rank
the schools, and prediction intervals can be constructed by resampling
from appropriate binomial distributions. The results in Figure
\ref{fig5} indicate the increased confidence we can have in the upper
extreme, with the top school identified with reasonable certainty. In
this example, the possible range of scores for ranking has finite
support, being restricted to the interval $[0,1]$; thus it is a context
where Theorem \ref{theo2} is applicable.\looseness=1

%
\begin{figure}

\includegraphics{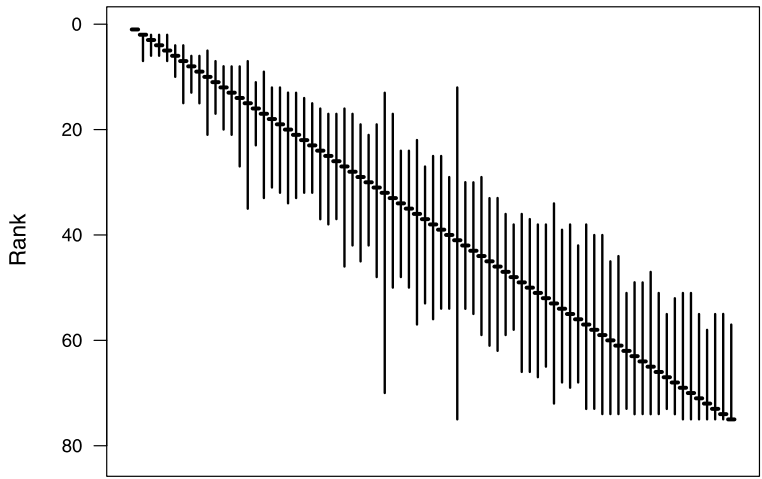}

\caption{Rankings of schools by students' exam performance with
prediction intervals.}
\label{fig5}
\end{figure}

Hill's (\citeyear{H75}) estimator of $\alpha$, when (\ref{eq2.2})
holds, is
relatively stable in this example and suggests that $\alpha\approx6$.
From (\ref{eq2.9}), we can calculate that $(n^{\alpha
/2}/p)^{1/(2\alpha-1)}
\approx4$, which is consistent with a small number of schools being
correctly ranked. If the number were large, then we would expect a
significant portion of the schools to be ranked with a high degree of
accuracy. In the case of these data, however, the small value suggests
that it might not be possible to obtain any correct ranks.
\end{Example}
\begin{Example}[(Simulation with exponential tails and infinite
support)]\label{Example4}
Here, we simulate increasing $n$ and $p$ in the case of
exponential tails. For a given $n$, set $p=0.0005 n^2$, let the $\Theta
_j$'s be drawn from a standard exponential distribution and the
$Q_{ij}$'s be normal random variables with zero mean and standard
deviation~3.5. Table \ref{tab1} shows the results of 1000 simulations
for various values of $n$, approximating (\ref{eq1.5}) for different
%
%
\begin{table}[b]
\caption{Probability that the first $j_0$ rankings are correct
in the case of exponential tails} \label{tab1}
\begin{tabular*}{\tablewidth}{@{\extracolsep{\fill
}}lcd{1.4}ccd{1.4}d{1.4}d{1.4}@{}}
\hline
& \multicolumn{7}{c@{}}{$\bolds n$}\\[-4pt]
& \multicolumn{7}{c@{}}{\hrulefill}\\
$\bolds{j_0}$ & \multicolumn{1}{c}{\textbf{500}}
& \multicolumn{1}{c}{\textbf{1000}} & \multicolumn{1}{c}{\textbf{2000}}
& \multicolumn{1}{c}{\textbf{5000}} & \multicolumn{1}{c}{\textbf{10,000}}
& \multicolumn{1}{c}{\textbf{20,000}} & \multicolumn{1}{c@{}}{\textbf
{50,000}} \\
\hline
$1$ &0.909&0.9365&0.959&0.970&0.9745&0.9840&0.9910\\
$n^{0.15}$&0.764&0.823&0.767&0.844&0.897&0.872&0.890\\
$n^{0.20}$&0.591&0.700&0.655&0.683&0.667&0.664&0.743\\
$n^{0.25}$&0.420&0.406&0.424&0.383&0.334&0.402&0.428\\
$n^{0.30}$&0.183&0.188&0.180&0.116&0.101&0.079&0.069\\
$n^{0.35}$&0.056&0.030&0.021&0.004&0.002&0.000&0.001\\
\hline
\end{tabular*}
\end{table}
choices of $j_0$. Theorem \ref{theo1} suggests that the results should converge
to 1 if $j_0 = o(n^{1/4})$, and degrade, otherwise. This appears
consistent with the results. The difficulty of the problem due to the
quadratic growth of $p$ and the large error in $Q_{ij}$ is also
evident; even when $j_0=1$ and $n$ is large, reliable prediction of the
top rank is not assured.
\end{Example}
\begin{Example}[(Simulation with polynomial tails and infinite
support)]\label{Example5}
We use the same setup as in the previous example, except that
the generating distribution for the $\Theta_j$'s is Pareto,
$F(x)=1-x^{-\alpha
}$ for $x\geq1$, with $\alpha=4$. Theorem \ref{theo1} and (\ref
{eq2.4}) suggest
that the rate $n^{4/18} p^{1/9} = n^{4/9}$ is critical for $j_0$, and
this is consistent with the results in Table \ref{tab2}. This is an
%
%
\begin{table}
\caption{Probability that the first $j_0$ rankings are correct
in the case of exponential tails} \label{tab2}
\begin{tabular*}{\tablewidth}{@{\extracolsep{\fill}}lccccccc@{}}
\hline
& \multicolumn{7}{c@{}}{$\bolds{n}$}\\ [-4pt]
& \multicolumn{7}{c@{}}{\hrulefill} \\
$\bolds{j_0}$ & \multicolumn{1}{c}{\textbf{500}} & \multicolumn
{1}{c}{\textbf{1000}}
& \multicolumn{1}{c}{\textbf{2000}} & \multicolumn{1}{c}{\textbf{5000}}
& \multicolumn{1}{c}{\textbf{10,000}} & \multicolumn{1}{c}{\textbf{20,000}}
& \multicolumn{1}{c@{}}{\textbf{50,000}} \\
\hline
$(1/5) n^{0.35}$&0.884&0.832&0.908&0.920&0.898&0.921&0.945\\
$(1/5) n^{0.40}$&0.694&0.672&0.708&0.731&0.801&0.786&0.803\\
$(1/5) n^{4/9}$ &0.477&0.510&0.586&0.568&0.569&0.520&0.540\\
$(1/5) n^{0.50}$&0.283&0.242&0.252&0.161&0.140&0.120&0.096\\
$(1/5) n^{0.55}$&0.071&0.086&0.031&0.020&0.006&0.002&0.001\\
\hline
\end{tabular*}
\end{table}
%
%
\begin{table}[b]
\caption{Probability all ranks identified correctly when $\Theta
_j$ is uniformly distributed}\label{tab3}
\begin{tabular*}{\tablewidth}{@{\extracolsep{\fill}}lccccccc@{}}
\hline
& \multicolumn{7}{c@{}}{$\bolds{n}$}\\ [-4pt]
& \multicolumn{7}{c@{}}{\hrulefill} \\
$\bolds{k}$ & \multicolumn{1}{c}{\textbf{500}} & \multicolumn
{1}{c}{\textbf{1000}}
& \multicolumn{1}{c}{\textbf{2000}} & \multicolumn{1}{c}{\textbf{5000}}
& \multicolumn{1}{c}{\textbf{10,000}} & \multicolumn{1}{c}{\textbf{20,000}}
& \multicolumn{1}{c@{}}{\textbf{50,000}} \\
\hline
$1/6$&0.502&0.494&0.525&0.593&0.635&0.658&0.701\\
$1/5$&0.498&0.511&0.471&0.558&0.568&0.578&0.606\\
$1/4$&0.497&0.478&0.492&0.505&0.517&0.496&0.502\\
$1/3$&0.500&0.457&0.395&0.343&0.289&0.259&0.212\\
$1/2$&0.502&0.369&0.249&0.107&0.046&0.011&0.000\\
\hline
\end{tabular*}
\end{table}
easier problem than that in the previous example, because of the
polynomial decay of the tail. For instance, the top right-hand result
in the table suggests that the top nine ranks can be correctly
ascertained more than 90\% of the time when $p>50\mbox{,}000$, whereas the
figure 0.890 in the last column of Table \ref{tab1} suggests that, for
the distribution represented there, only the top five ranks have this
level of reliability.
\end{Example}
\begin{Example}[(Simulation with polynomial tails with finite
support)]\label{Example6}
Theorem~\ref{theo2} has many interesting
consequences, but the present example focuses on case~(iii), where
$\alpha>{\frac{1}{2}}$. First, let the $\Theta _j$'s be uniformly
distributed on $[0,1]$, and consider a case where the entire ranking is
correct. Using the notation of Section \ref{sec3} and taking
$\alpha=1$, Theorem \ref{theo2} implies that $p\asymp n^{1/4}$ defines
the critical growth in dimension. For simulation, we took $p=2 n^k$ for
various $k$, and scaled the (normally distributed) error for each $k$
such that the $n=500$ case had probability approximately $0.5$ of
correctly identifying all ranks. Each simulation was repeated 10,000
times, with results summarized in Table \ref{tab3}. As predicted,
growth rates in dimension slower than $n^{1/4}$ have probability of
correct ranking tending to 1, while those faster than $n^{1/4}$
degrade.

%
\begin{table}
\caption{Probability that lowest $10n^k$ scores identified
correctly}\label{tab4}
\begin{tabular*}{\tablewidth}{@{\extracolsep{\fill}}l@{ \ \ }c@{ \ \ }c@{ \ \ }c@{ \ \ }c@{ \ \ }c@{ \ \ }c@{ \ \ }c@{ \ \ }c@{}}
\hline
& \multicolumn{8}{c@{}}{$\bolds{n}$}\\ [-4pt]
& \multicolumn{8}{c@{}}{\hrulefill} \\
$\bolds{k}$ & $\bolds{5\times10^3}$ & $\bolds{1\times10^4}$
& $\bolds{2\times10^4}$ & $\bolds{5\times10^4}$ &
$\bolds{1\times10^5}$ & $\bolds{2\times10^5}$ & $\bolds{5\times10^5}$
& $\bolds{1\times10^6}$ \\
\hline
$0.05$&0.500&0.539&0.553&0.583&0.603&0.609&0.628&0.641\\
$0.07$&0.502&0.532&0.506&0.546&0.558&0.580&0.555&0.591\\
$1/11$&0.497&0.486&0.489&0.516&0.489&0.463&0.513&0.496\\
$0.11$&0.497&0.481&0.471&0.432&0.461&0.447&0.452&0.421\\
$0.13$&0.506&0.492&0.461&0.481&0.445&0.427&0.387&0.385\\
\hline
\end{tabular*}
\end{table}

Next, we examine the case $p=5\times10^{-6} n^2$, where dimension
grows at a quadratic rate; and $F(x) = x^\alpha$ on $[0,1]$, with
$\alpha=6$,
implying a reasonably severe tail. Theorem \ref{theo2} suggests that if
$j_0 =
o(p^{1/22})$, or equivalently if $j_0=o(n^{1/11})$, then (\ref{eq1.5})
should hold. Table \ref{tab4} shows the probability of ranking the
smallest $j_0 = 10 n^k$ scores correctly for various $k$ and $n$, with
10,000 simulations. Again the normal error is tuned so that the $n=5
000$ case has probability of close to ${\frac{1}{2}}$. The results
suggest that
$n^{1/11}$ indeed separates values of $k$ for which correct ranking is possible.
\end{Example}


\section{Technical arguments}\label{sec5}

We begin by giving a brief sketch of the proof of Theorem \ref{theo1}.
Two steps
in the proof are initially presented as lemmas, the first using
moderate deviation properties to approximate sums related to the object
of interest, and the second employing Taylor's expansion applied to R\'
{e}nyi representations of order statistics to show that the gaps
$\Theta
_{(j+1)} - \Theta_{(j)}$ have a high probability of being of reasonable
size. In the proof itself, we use Lemma \ref{lemma1} to bound the
probability in
(\ref{eq1.5}) from below [see (\ref{eqres1})] and then show that the
last two terms in this expression converge to zero, implying that the
probability converges to 1 if (\ref{eq2.5}) holds. For the converse,
assuming independence, we find an upper bound to the probability in
(\ref{eqres2}) and show that if this probability tends to one then the
sum~$s(n)$, introduced at (\ref{eqP.29}), must converge to zero, which
in turn implies (\ref{eq2.5}). Only the exponential tail case is
presented in detail; comments at the end of the proof describe the main
differences in the polynomial tail case.

Throughout, we let $\mathcal{E}(j_0)$ denote the event that
${\bar Q}_{R_j}+\Theta_{R_j}>{\bar Q}_{R_{j_0}}+\Theta_{R_{j_0}}$
for $j_0+1\leq
j\leq
p$, we define $\mathcal{E}_j$ to be the event that
$\Theta_{(j+1)}-\Theta_{(j)}\geq-({\bar Q}_{R_{j+1}}-{\bar
Q}_{R_j})$, and we take ${\tilde\mathcal{E}}(j_0)$
and ${\tilde\mathcal{E}}_j$ to be the respective complements. Also,
we let
$\zeta_j=\Theta_{(j+1)}-\Theta_{(j)}$ denote the $j$th gap, where
$\Theta_{(0)}=-\infty$ for convenience.


In Lemma \ref{lemma1} below, we write $\mathcal{O}$ to denote the
sigma-field generated
by the $\Theta_j$'s, $N$ for a standard normal random variable independent
of $\mathcal{O}$, $\delta_n$ for any given sequence of positive constants
$\delta_n$ converging to zero, and $\Delta$ for a generic random
variable satisfying $P(|\Delta|\leq\delta_n)=1$.
\begin{lemma}\label{lemma1} For any positive integer $j_0<p$, let
$\mathcal{J}$
denote the
set of positive, even integers less than or equal to $j_0$. Put
\[
T_{1j}={\min(\zeta_{j-1},\zeta_j)
\over2 (\operatorname{var} {\bar Q}_{R_j})^{1/2}} ,\qquad
T_{2j}={\zeta_j
\over\{\operatorname{var} ({\bar Q}_{R_{j+1}}-{\bar Q}_{R_j})\}
^{1/2}} .
\]
Then
%
%
\begin{eqnarray} \label{eqP.9}
&&\sum_{j=1}^{j_0}
P \biggl\{ |{\bar Q}_{R_j} |> {\frac{1}{2}}\min(\zeta_{j-1},\zeta_j ) \biggr\}
\nonumber\\[-8pt]\\[-8pt]
&&\qquad =2 \{1+o(1)\} \sum_{j=1}^{j_0}
P(|N|>T_{1j})+o(1).\nonumber
\end{eqnarray}
If in addition the components of the $Q_i$'s are independent then
%
%
\begin{eqnarray}\label{eqP.11}
&&E \biggl[\exp\biggl\{-\sum_{j\in\mathcal{J}} P ({\tilde\mathcal{E}}_j\mid
\mathcal{O}) \biggr\} \biggr]
\nonumber\\[-8pt]\\[-8pt]
&&\qquad \leq\{1+o(1)\} E \biggl[\exp\biggl\{-(1+\Delta)
\sum_{j\in\mathcal{J}} P(N>T_{2j}\mid\mathcal{O}) \biggr\} \biggr] .\nonumber
\end{eqnarray}
\end{lemma}
\begin{pf}
Using the arguments of \citet{RS65} and
\citet{A72}, it can be shown that, if the constant $C_2$ in (\ref
{eq2.3}) satisfies $C_2>B^2+2$ where $B>0$, then as $n$ (and hence,
also $p$) diverges,
%
%
\begin{eqnarray}
\label{eqP.5}
P \{|{\bar Q}_j|>x (\operatorname{var} {\bar Q}_j)^{1/2}\}
&=&\{1+o(1)\} 2 \{1-\Phi(x)\} , \\
\label{eqP.6}
\qquad P [- ({\bar Q}_{j_1}-{\bar Q}_{j_2} )
\geq x \{\operatorname{var} ({\bar Q}_{j_1}- {\bar Q}_{j_2} ) \}^{1/2}]
&=&\{1+o(1)\} \{1-\Phi(x)\} ,
\end{eqnarray}
uniformly in $0<x<B (\log p)^{1/2}$ and $j,j_1,j_2\geq1$ such that
$j_1\neq j_2$. Expression (\ref{eqP.6}) requires the independence
assumption. Therefore, since $C_2>2 (C_1+1)$ in~(\ref{eq2.3}), we can
take $B=(2+\varepsilon)^{1/2}$ for some $\varepsilon>0$, and then
(\ref{eqP.5}) and
(\ref{eqP.6}) hold uniformly in $0<x<\{(2+\varepsilon) \log p\}^{1/2}$.
Thus, as $n\to\infty$, they hold uniformly in all $x>0$, modulo an
$o(p^{-1})$ term. We use (\ref{eqP.5}) to derive (\ref{eqP.9}), while
(\ref{eqP.6}) implies that
\[
\sum_{j\in\mathcal{J}} P ({\tilde\mathcal{E}}_j )
= \{1+o(1)\} \sum_{j\in\mathcal{J}} P(N>T_{2j})+o(1) ,
\]
which leads to (\ref{eqP.11}).
\end{pf}
%
\begin{lemma}\label{lemma2} If (\ref{eq2.1}), indicating the case of exponential
tails, holds then there exist $B_4,B_5>0$ such that, for any choice of
constants $c_1,c_2$ satisfying $0<c_1<c_2<(4-\varepsilon)^{-1}$ with
$\varepsilon$ as
in (\ref{eq2.4b}), and for all $B_6>0$,
%
%
\begin{eqnarray}
\label{eqP.25}
\inf_{j\in[1,n^{c_1}]}
P \bigl\{\zeta_j Z_{j+1}^{-1}
(\log n)^{1-(1/\alpha)}\geq B_4 n^{-c_1} \bigr\} &=& 1-O (n^{-B_6} ) ,
\\
\label{eqP.26}
\inf_{j\in[n^{c_1},n^{c_2}]}
P \bigl\{B_4\leq j \zeta_j Z_{j+1}^{-1}
(\log n)^{1-(1/\alpha)} \leq B_5 \bigr\} &=& 1-O (n^{-B_6} ) .
\end{eqnarray}
Note further that the constraint on $c_2$ permits $n^{c_2}$ to be of
size $\nu_{\exp}n^{\varepsilon_1}$ (where $\varepsilon_1>0$).
\end{lemma}
\begin{pf}
If $U_{(1)}\leq\cdots\leq U_{(p)}$ denote the order statistics of a
sample of size $p$ drawn from the uniform distribution on $[0,1]$ then,
for each $p$, we can construct a collection of independent random
variables $Z_1,\ldots,Z_p$ with the standard negative exponential
distribution on $[0,\infty]$, such that, for $1\leq j\leq p$, $U_{(j)}
=1-\exp(-V_j)$ where
\[
V_j=\sum_{k=1}^j {Z_k\over p-k+1}=w_j+W_j .
\]
For details, see \citet{R53}. Further, uniformly in $1\leq j\leq
{\frac{1}{2}}p$ and $2\leq p<\infty$,
%
%
\begin{eqnarray}\label{eqP.13}\qquad
&\displaystyle w_j=\sum_{k=p-j+1}^p {1\over k}
={j\over p}+O (j^2 /p^2 )=O(j/p) ,&
\\
%
%
\label{eqP.14}
&\displaystyle W_j=\sum_{k=p-j+1}^p k^{-1}(Z_{p-k+1}-1) ,\qquad
\sup_{1\leq j\leq p/2} j^{-1/2}|W_j|\leq p^{-1}W(p) ,
&
\\
%
%
\label{eqP.15}
&\displaystyle
\sup_{1\leq j\leq p/2} j^{-3/2} \Biggl|W_j-{1\over p} \sum_{k=p-j+1}^p
(Z_{p-k+1}-1) \Biggr|\leq p^{-2}W(p) ,&
\end{eqnarray}
where the nonnegative random variable $W(p)$, which without loss of
generality, we take to be common to (\ref{eqP.14}) and (\ref
{eqP.15}), satisfies the expression $P\{W(p)>p^\varepsilon\}
=O(p^{-C})$ for
each $C,\varepsilon>0$.

Using the second identity in (\ref{eqP.13}), and (\ref{eqP.14}), we
deduce that
%
%
\begin{eqnarray} \label{eqP.16}\quad
U_{(j+1)}-U_{(j)}&=&(V_{j+1}-V_j) \biggl\{1-{\frac{1}{2}}(V_{j+1}+V_j)
\nonumber\\
&&\hspace*{60.3pt}{}
+{\frac{1}{6}} (V_{j+1}^2+V_j V_{j+1}+V_j^2)-\cdots\biggr\} \\
&=&{Z_{j+1}\over p-j}
\biggl\{1+\Psi_{j1} \biggl({j\over p}+{S_{j1}\over p^{1/2}} \biggr) \biggr\},\nonumber
\end{eqnarray}
uniformly in $1\leq j\leq{\frac{1}{2}}p$, where the random variable
$\Psi
_{j1}$ satisfies, for $k=1$,
%
%
\begin{equation}\label{eqP.17}
P \Bigl({\max_{1\leq j\leq p/2}} |\Psi_{jk}|\leq A \Bigr)=1 ,
\end{equation}
$A>0$ is an absolute constant, and for each $C,\varepsilon>0$ the nonnegative
random variable $S_{j1}$ satisfies, with $k=1$,
%
%
\begin{equation}\label{eqP.18}
P \Bigl(\sup_{1\leq j\leq p/2} S_{jk}>p^\varepsilon\Bigr)
=O (p^{-C} ) .
\end{equation}
Using the third identity in (\ref{eqP.13}) and (\ref{eqP.15}), we
deduce that
%
%
\begin{eqnarray} \label{eqP.19}
0&\leq
&U_{(j)}=w_j+W_j-{\frac{1}{2}}(w_j+W_j)^2+\cdots\nonumber\\[-8pt]\\[-8pt]
&=&{j\over p}+\Psi_{j2} \biggl({j^2\over p^2}
+{j^{1/2}S_{j2}\over p} \biggr) ,\nonumber
\end{eqnarray}
where $\Psi_{j2}$ and $S_{j2}\geq0$ satisfy (\ref{eqP.17}) and (\ref
{eqP.18}), respectively.

Define $D_j=U_{(j+1)}-U_{(j)}$ and without loss of generality, $C_0=1$ in
(\ref{eq2.1}). If the common distribution function of the $\Theta
_j$'s is
$F$, then by Taylor's expansion,
%
%
\begin{eqnarray}\label{eqP.22}
\zeta_j &=& F^{-1}\bigl(U_{(j)}+D_j\bigr)-F^{-1}\bigl(U_{(j)}\bigr) \nonumber\\
&=& D_j (F^{-1})' \bigl(U_{(j)}+\omega_j D_j\bigr) \\
&=& \Psi_j {D_j\over U_{(j)}+\omega_j D_j}
\bigl\{-\log\bigl(U_{(j)}+\omega_j D_j\bigr) \bigr\}^{(1/\alpha)-1} ,\nonumber
\end{eqnarray}
where $0\leq\omega_j\leq1$ and the last line makes use of (\ref
{eq2.1}). The random variable $\Psi_j$ satisfies, for constants $B_1$,
$B_2$ and $B_3$ satisfying $0<B_1<B_2<\infty$ and \mbox{$0<B_3<1$},
\[
P \bigl(B_1\leq\Psi_j\leq B_2 \mbox{ for all $j$ such that }
U_{(j+1)}<B_3 \bigr)=1 .
\]
The required result then follows from (\ref{eqP.16}), (\ref{eqP.19})
and (\ref{eqP.22}).
\end{pf}
%
\begin{pf*}{Proof of Theorem \ref{theo1}}
Take $j_0<p$ a positive integer. Note that, taking $\mathcal{E}(j_o)$,
$\mathcal{E}
_j$, ${\tilde\mathcal{E}}(j_o)$, ${\tilde\mathcal{E}}_j$, $\mathcal
{O}$ and $\mathcal{J}$ as for Lemma \ref{lemma1},
\begin{eqnarray*}
&&\{{\widehat R}_j=R_j \mbox{ for } 1\leq j\leq j_0 \}  \\
&&\qquad \supseteq\bigl\{ |{\bar Q}_{R_j} |\leq{\tfrac{1}{2}}\min(\zeta
_{j-1},\zeta_j )
\mbox{ for } 1\leq j\leq j_0 \bigr\}
\cap\mathcal{E}(j_0) ,
\end{eqnarray*}
where we define $\Theta_{(j-1)}=-\infty$ if $j=1$ as before. Therefore,
defining $\pi(j_0)=P({\widehat R}_j=R_j$ for $1\leq j\leq
j_0)$, we
deduce that
%
%
\begin{equation}\label{eqP.1}
\pi(j_0) \geq1-\sum_{j=1}^{j_0}
P \biggl\{ |{\bar Q}_{R_j} |>{\frac{1}{2}}\min(\zeta_{j-1},\zeta_j ) \biggr\}
-P\{{\tilde\mathcal{E}}
(j_0)\} .
\end{equation}

Also,
\begin{eqnarray*}
&&\{{\widehat R}_j=R_j \mbox{ for } 1\leq j\leq j_0 \} \\
&&\qquad= \bigl\{{\bar X}_{R_1}\leq\cdots\leq{\bar X}_{R_{j_0}} \mbox{ and }
{\bar X}_j>{\bar X}_{R_{j_0}} \mbox{ for } j\notin\{R_1,\ldots
,R_{j_0}\}
\bigr\} \\
&&\qquad= \bigl\{\zeta_j\geq- ({\bar Q}_{R_{j+1}}-{\bar Q}_{R_j} )
\mbox{ for } 1\leq j\leq j_0 \\
&&\hspace*{38pt}
\mbox{and } \Theta_j-\Theta_{(j_0)}\geq- ({\bar Q}_j-{\bar
Q}_{R_{j_0}} )
\mbox{ for } j\notin\{R_1,\ldots,R_{j_0}\} \bigr\},
\end{eqnarray*}
and so
%
%
\begin{equation}\label{eqP.2}
\pi(j_0)\leq P \{\zeta_j \geq- ({\bar Q}_{R_{j+1}}
-{\bar Q}_{R_j} ) \mbox{ for } 1\leq j\leq j_0 \} .
\end{equation}

Letting $\pi_1(j_0)$ denote the probability that $\mathcal{E}_j$
holds for
all $j\in\mathcal{J}$, by (\ref{eqP.2}),
%
%
\begin{equation}\label{eqP.3}
\pi(j_0)\leq\pi_1(j_0) .
\end{equation}

Note that if the components of each $Q_i$ are independent, then the
events $\mathcal{E}_j$, for $j\in\mathcal{J}$, are independent
conditional on $\mathcal{O}
$. Therefore,
%
%
\begin{eqnarray}\label{eqP.4}
\pi_1(j_0)
&=&E \biggl\{P \biggl(\bigcap_{j\in\mathcal{J}} \mathcal{E}_j\mid\mathcal{O}\biggr)
\biggr\}
=E \biggl[\prod_{j\in\mathcal{J}} \{1-P (\mathcal{E}_j\mid\mathcal
{O}) \} \biggr] \nonumber\\[-8pt]\\[-8pt]
&\leq& E \biggl[\exp\biggl\{-\sum_{j\in\mathcal{J}} P ({\tilde\mathcal
{E}}_j\mid\mathcal{O}) \biggr\}\biggr].\nonumber
\end{eqnarray}

Using Lemma \ref{lemma1}, we have the following inequalities regarding
$\pi(j_0)$:
%
%
\begin{eqnarray}
\label{eqres1}
\pi(j_0) &\ge& 1- 2 \{1+o(1)\} \sum_{j=1}^{j_0}
P(|N|>T_{1j}) - P\{{\tilde\mathcal{E}}(j_0)\}+o(1),\\
\label{eqres2}
\pi(j_0) &\le& \{1+o(1)\} E \biggl[\exp\biggl\{-(1+\Delta)
\sum_{j\in\mathcal{J}} P(N>T_{2j}\mid\mathcal{O}) \biggr\}\biggr] .
\end{eqnarray}

To show that (\ref{eq2.5}) implies (\ref{eq1.5}), by (\ref{eqres1})
it is sufficient to show that $P\{{\tilde\mathcal{E}}(j_0)\}$ and
$\sum_{j=1}^{j_0}
P(|N|>T_{1j})$ are both $o(1)$, which we shall do in turn.

Define $\ell=(\log n)^{(1/\alpha)-1}$, let $N$ be a standard normal
random variable independent of $\mathcal{O}$, and let $Z$ be
independent of
$N$ and have the standard negative exponential distribution. Let $K_1$
be a positive constant. If $a_n$ is a sequence of positive numbers and
$f_n$ is a sequence of nonnegative functions, write $a_n\doteq f_n(K)$
to mean that, for constants $L_1,L_2>1$, either (a) $a_n\leq L_1
f_n(K)$ whenever $K\geq L_2$ and $n$ is sufficiently large, and
$a_n\geq L_1^{-1}f_n(K)$ whenever $K\leq L_2^{-1}$ and $n$ is
sufficiently large, or (b) $a_n\geq L_1^{-1}f_n(K)$ whenever $K\geq
L_2$ and $n$ is sufficiently large, and $a_n\leq L_1 f_n(K)$ whenever
$K\leq L_2^{-1}$ and $n$ is sufficiently large. Let $0<c_1<c_2<{\frac{1}{2}}$
and $c_1<{\frac{1}{4}}$, and let $j_0$ and $j_1$ denote integers satisfying
$|j_1-n^{c_1}|\leq1$, $j_1\leq j_0\leq n^{c_2}$ and $j_1/j_0\to0$.

When (\ref{eq2.1}) holds with $C_0=1$, Lemma \ref{lemma2} implies that,
for each
$B_6>0$ and letting $\gamma_j = n^{-1/2}j \ell^{-1}$,
%
%
\begin{eqnarray} \label{eqP.29}\hspace*{1pt}
s(n) &\equiv& \sum_{j=1}^{j_0} P\{|N|>K_1 n^{1/2}
(\zeta_j)\} \nonumber\\
&\doteq& O \{j_1 P (|N|>K_2Z\gamma_{j_1}^{-1})
+n^{-B_6} \}
+\sum_{j_1<j\leq j_0} P (|N|>KZ\gamma_j^{-1}) \nonumber\\
&\doteq& O \biggl\{j_1 \biggl(P (Z\leq\gamma_{j_1} )
+E \biggl[Z^{-1}\gamma_{j_1}
\exp\biggl\{-{\frac{1}{2}}(K Z\gamma_{j_1}^{-1})^2 \biggr\}
I (Z>\gamma_{j_1} ) \biggr] \biggr) \biggr\}\nonumber\\[-8pt]\\[-8pt]
&&{} +\sum_{j_1<j\leq j_0}
\biggl(P (Z\leq\gamma_j )
+E \biggl[Z^{-1}\gamma_{j}
\exp\biggl\{-{\frac{1}{2}}(K Z\gamma_{j}^{-1})^2 \biggr\}
I (Z>\gamma_{j} ) \biggr] \biggr) \nonumber\\
&\doteq& O \biggl\{j_1 \biggl(\gamma_{j_1}
+E \biggl[Z^{-1}\gamma_{j_1}
\exp\biggl\{-{\frac{1}{2}}(K Z\gamma_{j_1}^{-1})^2 \biggr\}
I (Z>\gamma_{j_1} ) \biggr] \biggr) \biggr\}\nonumber\\
&&{} +\sum_{j_1<j\leq j_0} \biggl(\gamma_j
+E \biggl[Z^{-1}\gamma_{j}
\exp\biggl\{-{\frac{1}{2}}(K Z\gamma_{j}^{-1})^2 \biggr\}
I (Z>\gamma_{j}) \biggr] \biggr) . \nonumber
\end{eqnarray}
Now,
\begin{eqnarray*}
&&
E \biggl[Z^{-1}\gamma_j
\exp\biggl\{-{\frac{1}{2}}(KZ\gamma_j^{-1})^2 \biggr\}
I (Z>\gamma_j ) \biggr] \\
&&\qquad=\int_{\gamma_j}^\infty z^{-1}\gamma_j
\exp\biggl\{-{\frac{1}{2}}(KZ\gamma_j^{-1})^2-z \biggr\} \,dz \\
&&\qquad=\gamma_j\int_1^\infty u^{-1}
\exp\biggl\{-{\frac{1}{2}}(K u)^2-\gamma_j u \biggr\} \,du
\asymp\gamma_j \\
&&\qquad= n^{-1/2}j \ell.
\end{eqnarray*}
(Here, we have used the fact that $j\leq j_0\leq n^{c_2}$ where
$c_2<{\frac{1}{2}}$.) Therefore,
%
%
\begin{eqnarray} \label{eqP.30}
s(n)&\asymp& j_1\cdot n^{-1/2}j_1 \ell^{-1}
+\sum_{j_1<j\leq j_0} n^{-1/2}j \ell^{-1}\nonumber\\
&\asymp& n^{-1/2}j_1^2 \ell^{-1}+n^{-1/2}j_0^2 \ell^{-1}\\
&\asymp& n^{-1/2}j_0^2 \ell^{-1}.\nonumber
\end{eqnarray}
(Here, we have used the fact that $j_1/j_0\to0$.)

The right-hand side of (\ref{eqP.30}) converges to zero if and only if
(\ref{eq2.5}) holds. Moreover, in view of the fact that
\[
P(|N|>T_{1j}) \leq P \biggl(|N|>{\zeta_{j-1} \over2 (\operatorname{var}
{\bar Q}_{R_j})^{1/2}
} \biggr)+P \biggl(|N|>{\zeta_j
\over2 (\operatorname{var} {\bar Q}_{R_j})^{1/2}} \biggr) ,
\]
and depending on the choice of $K_1$ in the definition of $s(n)$ at
(\ref{eqP.29}), $s(n)$ can be an upper bound to the series $\sum
_{j=1}^{j_0} P(|N|>T_{1j})$ on the right-hand side of (\ref{eqP.9}). Hence,
%
%
\begin{equation} \label{res1.1}
\sum_{j=1}^{j_0} P(|N|>T_{1j}) =o(1) .
\end{equation}
This deals with the second term on the right-hand side of (\ref{eqres1}).
Similarly, if $r\in[2,\infty)$ is a fixed integer, and if
$j_0=o(n^{1/4} \ell^{1/2})$, then
%
%
\begin{equation}\label{eqP.31}
s_1(n)\equiv\sum_{j=j_0+1}^{j_0+r-1}
P \{|N|>K_1 n^{1/2}(\zeta_j) \}=o(1) .
\end{equation}
Moreover, if $j_1$ denotes the integer part of $n^{c_2}-j_0$ then, for
constants $K_2$ and $K_3$ satisfying $K_1>K_2>K_3>0$, and for any $B>0$,
%
%
\begin{eqnarray} \label{eqP.32}
s_2(n)&\equiv&\sum_{j=j_0+r}^{j_0+j_1}
P \bigl\{|N|>K_1 n^{1/2}\bigl(\Theta_{(j+1)}-\Theta_{(j_0)}\bigr) \bigr\} \nonumber\\
&\leq&\sum_{j=r}^{j_1}
P \Biggl\{|N|>K_2 n^{1/2}\ell\sum_{k=1}^j (j_0+k)^{-1}Z_k \Biggr\}
+O (n^{-B} ) \nonumber\\[-8pt]\\[-8pt]
&\leq& j_1 P \{|N|>K_2 n^{1/4} \ell^{1/2}(Z_1+\cdots+Z_r) \}
+O (n^{-B} ) \nonumber\\
&=&O \{j_1 (n^{1/2}\ell^2 )^{-r} \},\nonumber
\end{eqnarray}
where we have assumed that $j_0=o(n^{1/4} \ell^{1/2})$ and also used
the fact that $Z_1+\cdots+Z_r$ has a $\operatorname{gamma}(r,1)$ distribution. If we
choose $r$ so large that $p n^{-r/2}=O(n^{-\varepsilon})$ for some
$\varepsilon>0$,
then we can deduce from (\ref{eqP.31}) and (\ref{eqP.32}) that
$s_1(n)+s_2(n)\to0$, and hence, by (\ref{eqP.26}), that
%
%
\begin{equation}\label{eqP.33}
\sum_{j=j_0+1}^{n^{c_2}}
P({\bar Q}_{R_j}+\Theta_{R_j}>{\bar Q}_{R_{j_0}}+\Theta_{R_{j_0}}
)\to0 .
\end{equation}

A more crude argument can be used to prove that if $r$ is so large that
$p^2 n^{-r/2}=O(n^{-\varepsilon})$ for some $\varepsilon>0$, and if
$j_0=o(n^{1/4}
\ell^{1/2})$, then
%
%
\begin{equation}\label{eqP.34}
\sum_{n^{c_2}<j\leq p}
P({\bar Q}_{R_j}+\Theta_{R_j}>{\bar Q}_{R_{j_0}}+\Theta_{R_{j_0}}
)\to0 .
\end{equation}
Together, (\ref{eqP.33}) and (\ref{eqP.34}) imply that if
$j_0=o(n^{1/4} \ell^{1/2})$ then
%
%
\begin{equation} \label{eqP.35}
P\{{\tilde\mathcal{E}}(j_0)\}\to0 .
\end{equation}

Thus, in light of (\ref{eqres1}), we see (\ref{res1.1}) and (\ref
{eqP.35}) imply that (\ref{eq2.5}) is sufficient for (\ref{eq1.5}).

We next show that (\ref{eq1.5}) implies (\ref{eq2.5}) in the
independent case. If (\ref{eq1.5}) holds, then by (\ref{eqres2}),
\[
\sum_{j\in\mathcal{J}} P(N>T_{2j}\mid\mathcal{O})\to0
\]
in probability. Therefore, by Lemma \ref{lemma2}, with $j_0$ and $j_1$
as above,
there exists $K_1>0$ such that
\[
\sum_{j_1<j\leq j_0} P\{|N|>n^{1/2}K_1 (\zeta_j)\mid\mathcal{O}\}
\to0
\]
in probability. (We can take the sum over all $j\in[j_1+1,j_0]$,
rather than just over even $j$, since (\ref{eqP.11}) holds for sums
over odd $j$ as well as over even $j$.) Hence, arguing as in the lines
below (\ref{eqP.29}), we deduce that for sufficiently large $K_2>0$,
%
%
\begin{equation}\label{eqP.36}
\sum_{j_1<j\leq j_0} f(Z_j/\delta_j)\to0
\end{equation}
in probability, where the random variables $Z_j$ are independent and
have a common exponential distribution, $\delta_j=n^{-1/2}j \ell^{-1}
$ and
\[
f(z)=z^{-1}\exp(-K_2 z^2) I(z>1) .
\]
We claim that this implies that the expected value of the left-hand
side of (\ref{eqP.36}) also converges to 0:
%
%
\begin{equation}\label{eqP.37}
\sum_{j_1<j\leq j_0} E\{f(Z_j/\delta_j)\}\to0
\end{equation}
or equivalently that $\sum_{j_1<j\leq j_0} \delta_j\to0$, and
thence [using the argument leading to (\ref{eqP.30})] that $s(n)\asymp
n^{-1/2}j_0^2 \ell^{-1}\to0$, which is equivalent to (\ref{eq2.5}).
Therefore, if we establish (\ref{eqP.37}) then we shall have proved
that (\ref{eq1.5}) implies (\ref{eq2.5}).

It remains to show that (\ref{eqP.36}) implies (\ref{eqP.37}). This
we do by contradiction. If (\ref{eqP.37}) fails then, along a
subsequence of values of $n$, the left-hand side of (\ref{eqP.37})
converges to a nonzero number. For notational simplicity, we shall make
the inessential assumptions that the number is finite and that the
subsequence involves all $n$, and we shall take $K_2=1$ in the
definition of $f$. In particular,
%
%
\begin{equation}\label{eqP.38}
t(n)\equiv\sum_{j_1<j\leq j_0} E\{f(Z_j/\delta_j)\}\to t(\infty),
\end{equation}
where $t(\infty)$ is bounded away from 0. Now, $t(n)=\{1+o(1)\} \mu
(1) \delta(n)$, where $\delta(n)=\sum_{j_1<j\leq j_0} \delta
_j$ and, for general $\lambda\geq1$, $\mu(\lambda)=\int
_{z>\lambda}z^{-1}
\exp(-z^2) \,dz$. Therefore,
%
%
\begin{equation}\label{eqP.39}
\delta(n)\to\delta(\infty)\equiv t(\infty)/\mu(1) .
\end{equation}
For each $\lambda>1$ the left-hand side of (\ref{eqP.36}) equals
$\Delta
_1+\Delta_2$, where, in view of (\ref{eqP.38}),
%
%
\begin{eqnarray} \label{eqP.40}
E(\Delta_2)&=&\sum_{j_1<j\leq j_0} E\{f(Z_j/\delta_j) I(Z_j>\lambda
\delta_j)\}\nonumber\\[-8pt]\\[-8pt]
&=&\{1+o(1)\} \mu(\lambda) \delta(n)\nonumber
\end{eqnarray}
and
\[
\Delta_1
=\sum_{j_1<j\leq j_0} f(Z_j/\delta_j) I(Z_j\leq\lambda\delta_j)
=\sum_{j_1<j\leq j_0} f(W_j) I_j
\]
with $W_j=Z_j/\delta_j$ and $I_j=I(\delta_j\leq Z_j\leq\lambda
\delta_j)$. However,
\[
\sum_{j_1<j\leq j_0} P(I_j=1)=\mu_1(\lambda) \delta(n)+o(1)
=\delta(\infty) \mu_1(\lambda)+o(1) ,
\]
where $\mu_1(\lambda)=\int_{1<z<\lambda}z^{-1}\exp(-z^2) \,dz$. Therefore,
in the limit as $n\to\infty$, $\Delta_1$ equals a sum, $S_\lambda
$ say, of $N$
independent random variables each having the distribution of $f(W)$,
where $W$ is uniformly distributed on $[1,\lambda]$, $N$ has a Poisson
distribution with mean $\delta(\infty) \mu_1(\lambda)$, and $N$ and
the summands are independent. The distribution of $S_\lambda$ is
stochastically monotone increasing, in the sense that $P(S_\lambda>s)$
increases with $\lambda$. On the other hand, since $\mu(\lambda
)\to0$ as
$\lambda\to\infty$ then, by (\ref{eqP.39}) and (\ref{eqP.40}),
\[
\lim_{\lambda\to\infty} \limsup_{n\to\infty} E(\Delta_2)=0 .
\]
Combining these results, we deduce that $\Delta_1+\Delta_2$, that
is, the
left-hand side of (\ref{eqP.36}), does not converge to zero in
probability. This contradicts (\ref{eqP.36}) and so establishes that
$t(\infty)$ must equal zero; that is, (\ref{eqP.37}) holds.\vspace*{8pt}

\textit{Comments on proving the polynomial case}: the proof for the
case of polynomial tails proceeds similarly. The main difference is
that in the proof of Lemma \ref{lemma2} we use (\ref{eq2.2}) instead of
(\ref
{eq2.1}), which forces a factor of $p^{-1/\alpha}$ into the results
of the
lemma, rather than $(\log n)^{1-(1/\alpha)}$. This in turn implies that
$s(n) \asymp n^{-1/2} j_0^{2+1/\alpha}p^{-1/\alpha}$, entailing that
convergence occurs if (and, in the case of independence, only if) $j_0
= o(\nu_{\mathrm{pol}})$, as required.
\end{pf*}

\printaddresses

\end{document}